\documentclass[a4paper,11pt]{article}
\usepackage{amsfonts}
\usepackage{amstext}
\usepackage{amsthm}
\usepackage{amsmath}
\usepackage{amssymb}
\usepackage{graphicx}
\usepackage{latexsym}
\usepackage[latin1]{inputenc}
\newcommand{\R}{\Bbb{R}}
\newcommand{\N}{\Bbb{N}}

\newtheorem{teor}{Theorem}[section]
\newtheorem{propo}{Proposition}[section]

\newtheorem{lema}{Lemma}[section]

\newtheorem{rem}{Remark}[section]

\newcommand{\n}{\noindent}

\newcommand {\fim}{\rule{0.5em}{0.5em}}

\begin{document}

\title{Fractional Elliptic Systems with Nonlinearities
of Arbitrary Growth
}

\author{\textbf{Edir Junior Ferreira Leite \footnote{\textit{E-mail addresses}:
edirjrleite@ufv.br (E.J.F. Leite)}}\\ {\small\it Departamento de Matem\'{a}tica,
Universidade Federal de Vi\c{c}osa,}\\ {\small\it CCE, 36570-000, Vi\c{c}osa, MG, Brazil}
}
\date{}{

\maketitle

\markboth{abstract}{abstract}
\addcontentsline{toc}{chapter}{abstract}

\hrule \vspace{0,2cm}

\n {\bf Abstract}

In this paper we discuss the existence, uniqueness and regularity of solutions of the following system of coupled semilinear Poisson equations on a smooth bounded domain $\Omega$ in $\R^n$:
\[
\left\{
\begin{array}{llll}
\mathcal{A}^s u= v^p & {\rm in} \ \ \Omega\\
\mathcal{A}^s v = f(u) & {\rm in} \ \ \Omega\\
u= v=0 & {\rm on} \ \ \partial\Omega
\end{array}
\right.
\]
where $s\in (0, 1)$ and $\mathcal{A}^s$ denote spectral fractional Laplace operators. We assume that $1< p<\frac{2s}{n-2s}$, and the function $f$ is superlinear and with no growth restriction (for example $f(r)=re^r$); thus the system has a nontrivial solution. Another important example is given by $f(r)=r^q$. In this case, we prove that such a system admits at least one positive solution for a certain set of the couple $(p,q)$ below the critical hyperbola
\[
\frac{1}{p + 1} + \frac{1}{q + 1} = \frac{n - 2s}{n}
\]
whenever $n > 2s$. For such weak solutions, we prove an $L^\infty$ estimate of Brezis-Kato type and derive the regularity property of the weak solutions.
\vspace{0.5cm}
\hrule\vspace{0.2cm}

\section{Introduction and main result}

This work is devoted to the study of existence and uniqueness of solutions for nonlocal elliptic systems on bounded domains which will be described henceforth.

The spectral fractional Laplace operator $\mathcal{A}^{s}$, is defined in terms of the Dirichlet spectra of the Laplace operator on $\Omega$. Roughly, if $(\varphi_k)$ denotes a $L^2$-orthonormal basis of eigenfunctions corresponding to eigenvalues $(\lambda_k)$ of the Laplace operator with zero Dirichlet boundary values on $\partial \Omega$, then the operator $\mathcal{A}^s$ is defined as $\mathcal{A}^{s} u = \sum_{k=1}^\infty c_k \lambda_k^s \varphi_k$, where $c_k$, $k \geq 1$, are the coefficients of the expansion $u = \sum_{k=1}^\infty c_k \varphi_k$.

A closely related to (but different from) the spectral fractional Laplace operator $\mathcal{A}^{s}$ is the restricted fractional Laplace operator $(-\Delta)^s$, see \cite{compare, rafaella}. This is defined as

\[
(-\Delta)^{s}u(x) = C(n,s)\, {\rm P.V.}\int\limits_{\R^{n}}\frac{u(x)-u(y)}{\vert x-y\vert^{n+2s}}\; dy\, ,
\]
\n for all $x \in \R^{n}$, where P.V. denotes the principal value of the first integral and

\[
C(n,s) = \left(\int\limits_{\R^{n}}\frac{1-\cos(\zeta_{1})}{\vert\zeta\vert^{n+2s}}\; d\zeta\right)^{-1}
\]

\n with $\zeta = (\zeta_1, \ldots, \zeta_n) \in \R^n$.

We remark that $(-\Delta)^{s}$ is a nonlocal operator on functions compactly supported in $\R^n$.

Factional Laplace operators arise naturally in several different areas such as Probability, Finance, Physics, Chemistry and Ecology, see \cite{A, bucur}. These operators have attracted special attention during the last decade. An extension for spectral fractional operator was devised by Cabr\'{e} and Tan \cite{CT} and Capella, D\'{a}vila, Dupaigne, and Sire \cite{CDDS} (see Br\"{a}ndle, Colorado, de Pablo, and S\'{a}nchez \cite{BrCPS} and Tan \cite{T} also). Thanks to these advances, the boundary fractional problem

\begin{equation} 
\left\{
\begin{array}{rrll}
\mathcal{A}^{s} u &=& u^p & {\rm in} \ \ \Omega\\
u &=& 0  & {\rm on} \ \ \partial\Omega
\end{array}
\right.
\end{equation}
has been widely studied on a smooth bounded open subset $\Omega$ of $\mathbb{R}^n$, $n\geq 2$, $s \in (0,1)$ and $p > 0$. Particularly, a priori bounds and existence of positive solutions for subcritical exponents ($p < \frac{n + 2s}{n - 2s}$) has been proved in \cite{BrCPS, CT, choi, CK, T} and nonexistence results has also been proved in \cite{BrCPS, tan1, T} for critical and supercritical exponents ($p \geq \frac{n + 2s}{n - 2s}$). The regularity result has been proved in \cite{CS, CDDS, T, yang}.

When $s = 1/2$, Cabr\'{e} and Tan \cite{CT} established the existence of positive solutions for equations having nonlinearities with the subcritical growth, their regularity, the symmetric property, and a priori estimates of the Gidas-Spruck type by employing a blow-up argument along with a Liouville type result for the square root of the Laplace operator in the half-space. Then \cite{T} has the analogue to $1/2<s<1$. Br\"{a}ndle, Colorado, de Pablo, and S\'{a}nchez \cite{BrCPS} dealt with a subcritical concave-convex problem. For $f(u)=u^q$ with the critical and supercritical exponents $q\geq\frac{n+2s}{n-2s}$, the nonexistence of solutions was proved in \cite{BaCPS, tan1, T} in which the authors devised and used the Pohozaev type identities. The Brezis-Nirenberg type problem was studied in \cite{tan1} for $s = 1/2$ and \cite{BaCPS} for $0<s<1$. The Lemma's Hopf and Maximum Principe was studied in \cite{T}.

The result of this paper are:

\begin{teor}\label{teo1} Suppose that $2\leq n<4s$, $0< p<\frac{2s}{n-2s}$, $f\in C(\mathbb{R})$, and set $F(r)=\int_0^rf(t)dt$. If there exist constants 
\[
\theta>\left\{
\begin{array}{llll}
2 & {\rm if} \ \ p> 1\\
1+\frac{1}{p} & {\rm if} \ \ p\leq1
\end{array},
\right.
\]
and $r_0\geq 0$ such that $\theta F(r)\leq f(r)r$ for all $\vert r\vert\geq r_0$ and 
\[
f(r)=\left\{
\begin{array}{llll}
o(r) & {\rm if} \ \ p> 1\\
o(r^{1/p}) & {\rm if} \ \ p\leq1
\end{array},
\right.
\]
for $r$ near $0$. Then the system
\begin{equation}\label{1}
\left\{
\begin{array}{llll}
\mathcal{A}^s u= v^p & {\rm in} \ \ \Omega\\
\mathcal{A}^s v = f(u) & {\rm in} \ \ \Omega\\
u= v=0 & {\rm on} \ \ \partial\Omega
\end{array},
\right.
\end{equation}
has a nontrivial solution.
\end{teor}

\begin{rem}
Note that if $2\leq n<4s$ then $n=2$ and $s\in\left(\frac{1}{2},1\right)$ or $n=3$ and $s\in\left(\frac{3}{4},1\right)$.
\end{rem}

\begin{teor} \label{teo2} Suppose that $n\geq 4s$, $0<p\leq 1$ and
\[
\left\{
\begin{array}{llll}
0<q<\frac{n+4s}{n-4s} & {\rm if} \ \ n>4s\\
0<q & {\rm if} \ \ n=4s
\end{array}.
\right.
\]
Then the system
\begin{equation}\label{2}
\left\{
\begin{array}{llll}
\mathcal{A}^s u= v^p & {\rm in} \ \ \Omega\\
\mathcal{A}^s v = u^q & {\rm in} \ \ \Omega\\
u= v=0 & {\rm on} \ \ \partial\Omega
\end{array},
\right.
\end{equation}
has a positive solution. Moreover, if $pq < 1$, then the problem (\ref{2}) admits a unique positive solution.
\end{teor}

\begin{rem}
Suppose that $n\geq 4s$, $0<q\leq 1$ and
\[
\left\{
\begin{array}{llll}
0<p<\frac{n+4s}{n-4s} & {\rm if} \ \ n>4s\\
0<p & {\rm if} \ \ n=4s
\end{array}.
\right.
\]
Clearly we have a result analogous to the above theorem.
\end{rem}

\begin{rem}
When $p, q > 1$, a priori bounds and existence of positive solutions of (\ref{2}) have been derived in \cite{choi} provided that $p, q$ satisfy
\begin{equation} \label{Hip}
\frac{1}{p+1} + \frac{1}{q+1} > \frac{n - 2s}{n}\, .
\end{equation}
\end{rem}

\begin{rem}
In the case when $n\leq 4s$, the above Theorems cover the remaining cases below the critical hyperbola and when $pq\neq 1$. In the case when $n>4s$, figure 1 exemplifies the region that the above theorem covers.
\end{rem}

\begin{figure}[ht]
\centering
\includegraphics[scale=1.2]{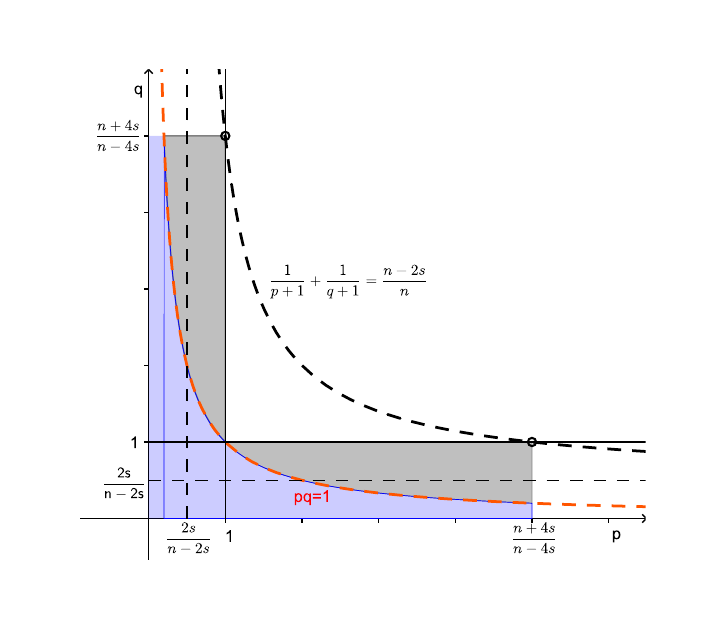}
\caption{The existence range of couples $(p,q)$ when $n>4s$.}
\end{figure}

\begin{rem}
For such weak solutions, we prove an $L^\infty$ estimate of Brezis-Kato type and derive the regularity property of the weak solutions based on the results obtained in \cite{CS} and \cite{T}.
\end{rem}

\begin{rem}
For $s = 1$, the problem (\ref{2}) and a number of its generalizations have been widely investigated in the literature, see for instance the survey \cite{DG} and references therein. Specifically, notions of sublinearity, superlinearity and criticality (subcriticality, supercriticality) have been introduced in \cite{FM, Mi1, Mi2, SZ1}. In fact, the behavior of (\ref{2}) is sublinear when $pq < 1$, superlinear when $pq > 1$ and critical (subcritical, supercritical) when $n \geq 3$ and $(p,q)$ is on (below, above) the hyperbola, known as critical hyperbola,

\[
\frac{1}{p+1}+\frac{1}{q+1}=\frac{n-2}{n}\, .
\]

\n When $pq = 1$, its behavior is resonant and the corresponding eigenvalue problem has been addressed in \cite{marcos}. The sublinear case has been studied in \cite{FM} where the existence and uniqueness of positive classical solution is proved. The superlinear-subcritical case has been covered in the works \cite{CFM}, \cite{DF}, \cite{DR} and \cite{vander} where the existence of at least one positive classical solution is derived. Lastly, the nonexistence of positive classical solutions has been established in \cite{Mi1} on star-shaped domains.
\end{rem}

\begin{rem}
When $0 < s < 1$ and $p, q > 0$, existence of positive solutions of (\ref{2}) for the restricted fractional Laplace operator $(-\Delta)^s$ have been derived in \cite{EM1} provided that $pq\neq 1$ and $(p, q)$ satisfies (\ref{Hip}).
\end{rem}

\begin{rem}
Related systems have been investigated by using other methods. We refer to the work \cite{EM} for systems involving different operators $(-\Delta)^{s}$ and $(-\Delta)^{t}$ in each one of equations. In this case, the study is non-variational.
\end{rem}

The rest of paper is organized into five sections. In Section 2 we briefly recall some definitions and facts related to fractional Sobolev spaces. In Section 3, we prove the case $p> 1$ of Theorem \ref{teo1} by applying the Strongly Indefinite Functional Theorem of Li-Willem. Then, we prove the case $p\leq 1$ by using the mountain pass theorem of Ambrosetti-Rabinowitz. In Section 4, we prove the case $pq<1$ of Theorem \ref{teo2} by using a direct minimization approach, Hopf lemma and maximum principles. Next we establish the remaining cases by using the mountain pass theorem. In Section 5, we establish regularity property of the weak solutions of system (\ref{1}) based on the results obtained in \cite{CS} and \cite{T}. Finally we establish the Brezis-Kato type result and derive the regularity of solutions to (\ref{2}).

\section{Preliminaries}

In this section we briefly recall some definitions and facts related to fractional Sobolev spaces.

We start by fixing a parameter $0 < s < 1$. Let $\Omega$ be an open subset of $\R^n$, with $n \geq 1$. For any $r \in (1,+\infty),$ one defines the fractional Sobolev space $W^{s,r}(\Omega)$ as

\begin{equation}\label{frac2.1}
W^{s,r}(\Omega):=\left\{u\in L^r(\Omega): \frac{\vert u(x)-u(y)\vert}{\vert x-y\vert^{\frac{n}{r}+s}}\in L^r(\Omega\times\Omega)\right\}\, ,
\end{equation}
that is, an intermediary Banach space between $L^r(\Omega)$ and $W^{1,r}(\Omega)$ induced with the norm
\begin{equation}\label{frac2.2}
\Vert u\Vert_{W^{s,r}(\Omega)}:=\left(\int\limits_{\Omega}\vert u\vert^r dx + \int\limits_{\Omega}\int\limits_{\Omega}\frac{\vert u(x)-u(y)\vert^r}{\vert x-y\vert^{n+sr}}dxdy\right)^{\frac{1}{r}}\, ,
\end{equation}
where the term
\[
[u]_{W^{s,r}(\Omega)}:=\left(\int\limits_{\Omega}\int\limits_{\Omega}\frac{\vert u(x)-u(y)\vert^r}{\vert x-y\vert^{n+sr}}dxdy\right)^{\frac{1}{r}}
\]
is the so-called Gagliardo semi-norm of $u$.

Let $s \in \R \setminus \N$ with $s\geq 1$. The space $W^{s,r}(\Omega)$ is defined as
\[
W^{s,r}(\Omega)=\{u\in W^{[s],r}(\Omega) : D^ju\in W^{s-[s],r}(\Omega),\forall j, \vert j\vert=[s]\}\, ,
\]

\n where $[s]$ is the largest integer smaller than $s$, $j$ denotes the $n$-uple $(j_1, \ldots, j_n) \in \N^n$ and $|j|$ denotes the sum $j_1 + \ldots + j_n$.

It is clear that $W^{s,r}(\Omega)$ endowed with the norm

\begin{equation}\label{frac2.11}
\Vert u\Vert_{W^{s,r}(\Omega)}=\left(\Vert u\Vert^r_{W^{[s],r}(\Omega)}+ [u]^r_{W^{s-[s],r}(\Omega)}\right)^{\frac{1}{r}}
\end{equation}
is a reflexive Banach space.

Clearly, if $s=m$ is an integer, the space $W^{s,r}(\Omega)$ coincides with the Sobolev space $W^{m,r}(\Omega)$.

Let $W^{s,r}_0(\Omega)$ denote the closure of $C^\infty_0(\Omega)$ with respect to the norm $\Vert\cdot\Vert_{W^{s,r}(\Omega)}$ defined in (\ref{frac2.11}). For $0 < s \leq 1$, we have
\[
W^{s,r}_0(\Omega)=\{u\in W^{s,r}(\R^n) : u=0 \text{ in }\R^n\setminus\Omega\}\, .
\]

In this paper, we focus on the case $r = 2$. This is quite an important case since the
fractional Sobolev spaces $W^{s,2}(\Omega)$ and $W^{s,2}_0(\Omega)$ turn out to be Hilbert spaces. They are usually denoted by $H^s(\Omega)$ and $H^s_0(\Omega)$, respectively.\\

The spectral fractional Laplace operator $\mathcal{A}^s$ is defined as follows. Let $\varphi_k$ be an eigenfunction of $-\Delta$ given by
\begin{equation}
\left\{\begin{array}{ccll}
-\Delta\varphi_k&=&\lambda_k\varphi_k & {\rm in} \ \ \Omega\; \\
\varphi_k&=& 0 & {\rm on} \ \ \partial\Omega
\end{array}\right.  ,
\end{equation}
where $\lambda_k$ is the corresponding eigenvalue of $\varphi_k,0<\lambda_1<\lambda_2\leq\lambda_3\leq\cdots\leq\lambda_k\rightarrow +\infty.$ Then, $\{\varphi_k\}_{k=1}^{\infty}$ is an orthonormal basic of $L^2(\Omega)$ satisfying 
\[
\int\limits_{\Omega}\varphi_j\varphi_kdx=\delta_{j,k}.
\]
We define the operator $\mathcal{A}^s$ for any $u\in C^\infty_0(\Omega)$ by

\begin{equation}\label{lapla}
\mathcal{A}^su=\sum_{k=1}^{\infty}\lambda_k^s\xi_k\varphi_k,
\end{equation}
where
\[
u=\sum_{k=1}^{\infty}\xi_k\varphi_k\text{  and  }\xi_k=\int\limits_{\Omega}u\varphi_kdx.
\]

This operator is defined on a Hilbert space
\[
\Theta^{s}(\Omega)=\{u=\sum_{k=1}^{\infty}u_k\varphi_{k}\in L^2(\Omega)\mid\sum_{k=1}^{\infty}\lambda^s_{k}\vert u_k\vert^2<+\infty\}
\]
with values in its dual $\Theta^{s}(\Omega)'$. Thus the inner product of $\Theta^{s}(\Omega)$ is given by
\[
\langle u,v\rangle_{\Theta^{s}(\Omega)}=\int\limits_{\Omega}\mathcal{A}^{s/2}u\mathcal{A}^{s/2}vdx=\int\limits_{\Omega}u\mathcal{A}^svdx=\int\limits_{\Omega}v\mathcal{A}^sudx.
\]
We denote by $\Vert\cdot\Vert_{\Theta^{s}}$ the norm derived from this inner product. We remark that $\Theta^{s}(\Omega)'$ can be described as the completion of the finite sums of the form
\[
f=\sum_{k=1}^{\infty}c_k\varphi_{k} 
\]
with respect to the dual norm
\[
\Vert f\Vert_{\Theta^{s}(\Omega)'}=\sum_{k=1}^{\infty}(\lambda_{k}^s)^{-1}\vert c_k\vert^2=\Vert\mathcal{A}^{-s/2}f\Vert_{L^2}^2=\int\limits_{\Omega}f\mathcal{A}^{-s}fdx
\]
and it is a space of distributions. Moreover, the operator $\mathcal{A}^s$ is an isomorphism between $\Theta^{s}(\Omega)$ and $\Theta^{s}(\Omega)'\simeq\Theta^{s}(\Omega)$, given by its action on the eigenfunctions. If $u,v\in\Theta^{s}(\Omega)$ and $f = \mathcal{A}^su$ we have, after this isomorphism,
\[
\langle f, v\rangle_{\Theta^{s}(\Omega)'\times\Theta^{s}(\Omega)} = \langle u, v\rangle_{\Theta^{s}(\Omega)\times\Theta^{s}(\Omega)} =\sum_{k=1}^{\infty}\lambda^s_{k}u_kv_k.
\]
If it also happens that $f\in L^{2}(\Omega)$, then clearly we get 
\[
\langle f, v\rangle_{\Theta^{s}(\Omega)'\times\Theta^{s}(\Omega)}=\int\limits_{\Omega}fv dx.
\]
We have $\mathcal{A}^{-s}:\Theta^{s}(\Omega)'\rightarrow\Theta^{s}(\Omega)$ can be written as
\[
\mathcal{A}^{-s}f(x)=\int\limits_{\Omega}G_{\Omega}(x,y)f(y)dy,
\]
where $G_{\Omega}$ is the Green function of operator $\mathcal{A}^s$ (see \cite{Green, ros176}). It is known that
\[
\Theta^s(\Omega)=\left\{
\begin{array}{llll}
L^2(\Omega) & {\rm if} \ \ s=0\\
H^{s}(\Omega)=H^{s}_0(\Omega) & {\rm if} \ \ s\in(0,\frac{1}{2})\\
H^{\frac{1}{2}}_{00}(\Omega) & {\rm if} \ \ s=\frac{1}{2}\\
H^{s}_0(\Omega) & {\rm if} \ \ s\in(\frac{1}{2},1]\\
H^{s}(\Omega)\cap H^{1}_0(\Omega) & {\rm if} \ \ s\in(1,2]
\end{array},
\right.
\]
where $H^{\frac{1}{2}}_{00}(\Omega):=\{u\in H^{1/2}(\Omega)\mid\int\limits_{\Omega}\frac{u^2(x)}{d(x)}dx<+\infty\}.$\\

Observe that the injection $\Theta^s(\Omega)\hookrightarrow H^{s}(\Omega)$ is continuous. By the Sobolev imbedding theorem we therefore have continuous imbeddings $\Theta^s(\Omega)\subset L^{p+1}(\Omega)$ if $p+1\leq\frac{2n}{n-2s}$ and these imbedding are compact if $p+1<\frac{2n}{n-2s}$ for $0<s<2n$. Also, we have compact imbedding  $\Theta^s(\Omega)\subset C(\Omega)$, if 
\[
\frac{s}{n}>\frac{1}{2}.
\]
For $0<r<2$ we have $\mathcal{A}^s:\Theta^{r}(\Omega)\rightarrow\Theta^{r-2s}(\Omega)$ is an isomorphism (see \cite{vander}).

Finally, by weak solutions, we mean the following: Let $f\in L^{\frac{2n}{n+2s}}(\Omega)$. Given the problem
\begin{equation}\label{prob2}
\left\{\begin{array}{ccll}
\mathcal{A}^{s}u&=& f & {\rm in} \ \ \Omega\; \\
u&=& 0 & {\rm on} \ \ \partial\Omega
\end{array}\right.  
\end{equation}
we say that a function $u\in\Theta^s(\Omega)$ is a weak solution of (\ref{prob2}) provided
\[
\int\limits_{\Omega}\mathcal{A}^{s/2}u\mathcal{A}^{s/2}\phi dx=\int\limits_{\Omega}f\phi dx
\]
for all $\phi\in\Theta^s(\Omega)$.

\section{Proof of Theorem \ref{teo1}}

We organize the proof of Theorem \ref{teo1} into two parts. We start by proving the existence of a
weak solution in case $p>1$. 

\subsection{The case $p>1$}

We define the product Hilbert spaces
\[
E^{\alpha}(\Omega)=\Theta^{\alpha}(\Omega)\times\Theta^{2s-\alpha}(\Omega),\text{   }0<\alpha<2s
\]
where your inner product is given by
\[
\langle(u_1,v_1),(u_2,v_2)\rangle_{E^{\alpha}(\Omega)}=\langle\mathcal{A}^{\alpha/2}u_1,\mathcal{A}^{\alpha/2}u_2\rangle_{L^{2}(\Omega)}+\langle\mathcal{A}^{s-\alpha/2}v_1,\mathcal{A}^{s-\alpha/2}v_2\rangle_{L^{2}(\Omega)}.
\]

We denote by $\Vert\cdot\Vert_{E}$ the norm derived from this inner product, i.e,

\[
\Vert(u,v)\Vert_{E}=\left(\Vert u\Vert^2_{\Theta^{\alpha}}+\Vert v\Vert^2_{\Theta^{2s-\alpha}}\right)^{\frac{1}{2}}.
\]

We also have $\mathcal{A}^s:\Theta^\alpha(\Omega)\rightarrow\Theta^{\alpha-2s}(\Omega)$ is an isomorphism, see \cite{vander}. Hence
$$\left(\begin{array}{ccllrr}
0& \mathcal{A}^{s}  \\
\mathcal{A}^{s} & 0 
\end{array}\right): E^\alpha(\Omega)\rightarrow\Theta^{-\alpha}\times\Theta^{\alpha-2s}(\Omega)=E^\alpha(\Omega)'$$
is an isometry. We consider the Lagrangian
\begin{eqnarray*}
\mathcal{J}(u,v)&=&\int\limits_{\Omega}\mathcal{A}^{s/2}u\mathcal{A}^{s/2}vdx-\int\limits_{\Omega}\left(\frac{1}{p+1}\vert v\vert^{p+1}+F(u)\right)dx,
\end{eqnarray*}
i.e., a strongly indefinite functional. The o Hamiltonian is given by
\begin{equation}\label{2.2}
\mathcal{H}(u,v)=\int\limits_{\Omega}\left(\frac{1}{p+1}\vert v\vert^{p+1}+F(u)\right)dx.
\end{equation}
The quadratic part can again be written as
\[
A(u,v)=\frac{1}{2}\langle L(u,v),(u,v)\rangle_{E^{\alpha}(\Omega)}=\int\limits_{\Omega}\mathcal{A}^{\alpha/2}u\mathcal{A}^{s-\alpha/2}vdx=\int\limits_{\Omega}\mathcal{A}^{s/2}u\mathcal{A}^{s/2}vdx,
\]
where
$$L=\left(\begin{array}{ccllrr}
0& \mathcal{A}^{s-\alpha}  \\
\mathcal{A}^{\alpha-s} & 0 
\end{array}\right)$$
is bounded and self-adjoint. Introducing the "diagonals"
\[
E^+=\{(u,\mathcal{A}^{\alpha-s}u):u\in\Theta^\alpha(\Omega)\}\text{ and }E^-=\{(u,-\mathcal{A}^{\alpha-s}u):u\in\Theta^\alpha(\Omega)\}
\]
we have
\[
E^{\alpha}(\Omega)=E^+\oplus E^-
\]

An orthonormal basis of $E^\alpha(\Omega)$ is given by
\[
\left\{\frac{1}{\sqrt{2}}(\lambda_k^{-\alpha/2}\varphi_k,\pm\lambda_k^{\alpha/2-s}\varphi_k):k=1,2,\cdots\right\}.
\]

The derivative of $A(u,v)$ defines a bilinear form
\begin{equation}\label{1.27}
B((u_1,v_1),(u_2,v_2))=A'(u_1,v_1)(u_2,v_2)=\langle L(u_1,v_1),(u_2,v_2)\rangle_{E^{\alpha}(\Omega)},
\end{equation}
where $(u_1,v_1),(u_2,v_2)\in E^{\alpha}(\Omega)$ with
\begin{equation}\label{1.28}
A(u_1,v_1)=\frac{1}{2}B((u_1,v_1),(u_1,v_1))\text{ and }B((u_1,v_1)^+,(u_1,v_1)^-)=0.
\end{equation}

We will give the choice of $\alpha$ in the following lemma.

\begin{lema}\label{lema3.1} Let $1 < p<\frac{2s}{n-2s}$. Then there exist parameter $0<\alpha<2s$ such that the following embeddings are continuous and compact:
\[
\Theta^{2s-\alpha}(\Omega)\subset L^{p+1}(\Omega) \text{  and  } \Theta^\alpha(\Omega)\subset C(\Omega).
\]
\end{lema}

{\bf Proof.} Note that $\Theta^{2s-\alpha}(\Omega)\subset L^{q}(\Omega)$ compactly, if $q<\frac{2n}{n-4s+2\alpha}$ and $\Theta^\alpha(\Omega)\subset C(\Omega)$ compactly, if 
\[
\frac{\alpha}{n}>\frac{1}{2}.
\]

We have $p+1<\frac{n}{n-2s}$. Thus if $\alpha>\frac{n}{2}$, then $p+1<\frac{2n}{n-4s+2\alpha}$.

For $n=2$, we have $s\in\left(\frac{1}{2},1\right)$. In this case choose $1<\alpha<2s$. For $n=3$, we have $s\in\left(\frac{3}{4},1\right)$. In this case choose $\frac{3}{2}<\alpha<2s$. This ends the proof. \; \fim

\begin{rem}
Note that $\alpha-s>0$. Thus $\Theta^\alpha(\Omega)\hookrightarrow\Theta^s(\Omega)$ is compact.
\end{rem}

The functional $\mathcal{J}(u, v): E^\alpha(\Omega)\rightarrow\mathbb{R}$ is strongly indefinite near zero, in the sense that there exist infinite dimensional subspaces $E^+$ and $E^-$ with $E^+\oplus E^- = E^\alpha(\Omega)$ such that the functional is (near zero) positive definite on $E^+$ and negative definite
on $E^-$. Li-Willem \cite{LI} prove the following general existence theorem for such situations, which can be applied in our case:

\begin{teor}\label{teo2.4} (Li-Willem, 1995). Let $\Phi:E\rightarrow\mathbb{R}$ be a strongly indefinite $C^1$-functional satisfying

\begin{itemize}
\item[(i)] $\Phi$ has a local linking at the origin, i.e. for some $r > 0$:
\[
\Phi(z)\geq 0 \text{ for }z\in E^+,\text{  } \Vert z\Vert_E\leq r \text{ and }\Phi(z)\leq 0,\text{ for }z\in E^-,\text{  } \Vert z\Vert_E\leq r;
\]
\item[(ii)] $\Phi$ maps bounded sets into bounded sets;
\item[(iii)] let $E^+_n$ be any $n$-dimensional subspace of $E^+$; then $\Phi(z)\rightarrow -\infty$ as $\Vert z\Vert\rightarrow +\infty$, $z\in E^+_n\oplus E^-$;
\item[(iv)] $\Phi$ satisfies the Palais-Smale condition $(PS)$ (Li-Willem \cite{LI} require a weaker $"(PS^\ast)$-condition", however, in our case the classical $(PS)$ condition will be satisfied).
\end{itemize}
 
Then $\Phi$ has a nontrivial critical point.
\end{teor}

We now verify that our functional satisfies the assumptions of this theorem.

First, it is clear, with the choice of $\alpha$ (Lemma \ref{lema3.1}), that $\mathcal{J}(u,v)$ is a $C^1$-functional on $E^\alpha(\Omega)$.

We show that the condition (i) of Theorem \ref{teo2.4} is satisfied. It is easy to see that $\mathcal{J}(u,v)$ has a local linking with respect to $E^+$ and $E^-$ at the origin.

Now the condition (ii) of Theorem \ref{teo2.4}. Let $B\subset E^\alpha(\Omega)$ be a bounded set, i.e. $\Vert u\Vert_{\Theta^\alpha}\leq c$, $\Vert v\Vert_{\Theta^{2s-\alpha}}\leq c$, for all $(u,v)\in B$.
Then
\begin{eqnarray*}
\vert\mathcal{J}(u,v)\vert &\leq & \Vert\mathcal{A}^{\alpha/2}u\Vert_{L^2}\Vert\mathcal{A}^{s-\alpha/2}v\Vert_{L^2}+\int\limits_{\Omega}\vert v\vert^{p+1}dx + \int\limits_{\Omega}\vert f(u)\vert dx\\
&\leq & \Vert u\Vert_{\Theta^\alpha}\Vert v\Vert_{\Theta^{2s-\alpha}}+c\Vert v\Vert^{p+1}_{\Theta^{2s-\alpha}} + \vert\Omega\vert\sup\{\vert f(u(x))\vert:x\in\Omega\}\leq C.
\end{eqnarray*}

We show that the condition (iii) of Theorem \ref{teo2.4} is satisfied. Let $z_k = z^+_k + z^-_k\in E^+_n\oplus E^-$ denote a sequence with $\Vert z_k\Vert_E\rightarrow +\infty$. By the above, $z_k$ may be written as
\[
z_k = (u_k,\mathcal{A}^{\alpha-s}u_k) + (w_k,-\mathcal{A}^{\alpha-s}w_k),
\]
with $u_k\in \Theta^\alpha_n(\Omega)$, $w_k\in \Theta^\alpha(\Omega)$,
where $\Theta^\alpha_n(\Omega)$ denotes an $n$-dimensional subspace of $\Theta^\alpha(\Omega)$. Thus, the functional $\mathcal{J}(z_k)$ takes the form
\begin{eqnarray*}
\mathcal{J}(z_k)&=&\int\limits_\Omega \mathcal{A}^{\alpha/2}u_k\mathcal{A}^{s-\alpha/2}\mathcal{A}^{\alpha-s}u_kdx -\int\limits_\Omega \mathcal{A}^{\alpha/2}w_k\mathcal{A}^{s-\alpha/2}\mathcal{A}^{\alpha-s}w_kdx \\
& &-\frac{1}{p+1}\int\limits_\Omega\vert\mathcal{A}^{\alpha -s}(u_k-w_k)\vert^{p+1}dx-\int\limits_\Omega F(u_k-w_k)dx\\
&=& \int\limits_\Omega \vert\mathcal{A}^{\alpha/2}u_k\vert^2dx -\int\limits_\Omega \vert\mathcal{A}^{\alpha/2}w_k\vert^2dx -\frac{1}{p+1}\int\limits_\Omega\vert\mathcal{A}^{\alpha -s}(u_k-w_k)\vert^{p+1}dx\\
& &-\int\limits_\Omega F(u_k-w_k)dx.
\end{eqnarray*}

Note that
\[
\Vert z_k\Vert_E\rightarrow\infty\Longleftrightarrow\int\limits_\Omega\vert\mathcal{A}^{\alpha/2}u_k\vert^2dx + \int\limits_\Omega\vert\mathcal{A}^{\alpha/2}w_k\vert^2dx =\Vert u_k\Vert^2_{\Theta^\alpha}+\Vert w_k\Vert^2_{\Theta^\alpha}\rightarrow\infty.
\]

Now, if

1) $\Vert u_k\Vert^2_{\Theta^\alpha}\leq c$, then $\Vert w_k\Vert^2_{\Theta^\alpha}\rightarrow\infty$, and then $\mathcal{J}(z_k)\rightarrow -\infty$;

2) $\Vert u_k\Vert^2_{\Theta^\alpha}\rightarrow +\infty$, then we estimate ($c, c_1$ and $c_2$ are positive constants) using the fact that $\alpha-s>0$ and $p > 1$
\[
\int\limits_\Omega\vert\mathcal{A}^{\alpha-s}(u_k - w_k)\vert^{p+1}dx\geq c\left(\int\limits_\Omega\vert\mathcal{A}^{\alpha-s}(u_k - w_k)\vert^2dx\right)^{\frac{p+1}{2}}\geq c_1\Vert u_k-w_k\Vert_{L^2}^{p+1}
\]
and
\[
\int\limits_\Omega F(u_k + w_k)dx\geq c_2\int\limits_\Omega\vert u_k + w_k\vert^{p+1}dx-d\geq c_1\Vert u_k+w_k\Vert_{L^2}^{p+1}-d
\]
and hence we obtain the estimate
\[
\mathcal{J}(z_k)\leq\frac{1}{2}\Vert u_k\Vert^2_{\Theta^\alpha}-c_1\left(\Vert u_k-w_k\Vert_{L^2}^{p+1}+\Vert u_k+w_k\Vert_{L^2}^{p+1}\right)+d.
\]
Since $\phi(t) = t^{p+1}$ is convex, we have $\frac{1}{2}(\phi(t)+\phi(r))\geq\phi\left(\frac{1}{2}(r+t)\right)$, and hence
\begin{eqnarray*}
\mathcal{J}(z_k) & \leq &\frac{1}{2}\Vert u_k\Vert^2_{\Theta^\alpha}-c_1\frac{1}{2^p}\left(\Vert u_k-w_k\Vert_{L^2}+\Vert u_k+w_k\Vert_{L^2}\right)^{p+1}+d \\
& \leq &\frac{1}{2}\Vert u_k\Vert^2_{\Theta^\alpha}-c_1\frac{1}{2^p}\Vert u_k\Vert^{p+1}_{L^2}+d.
\end{eqnarray*}
Since on $\Theta_n^\alpha(\Omega)$ the norms $\Vert u_k\Vert_{\Theta^\alpha}$ and $\Vert u_k\Vert_{L^2}$ are equivalent, we conclude that also in this case $\mathcal{J}(z_k)\rightarrow -\infty$.

Finally, the condition (iv) of Theorem \ref{teo2.4}. Let $(z_n)\subset E^\alpha(\Omega)$ denote a $(PS)$-sequence, i.e. such that
\begin{equation}\label{2.5}
\vert\mathcal{J}(z_n)\vert\rightarrow c,\text{ and }\vert(\mathcal{J}'(z_n),\eta)\vert\leq\varepsilon_n\Vert\eta\Vert_E,\forall\eta\in E^\alpha(\Omega),\text{ and }\varepsilon_n\rightarrow 0.
\end{equation}
We first show:

\begin{lema} 
The $(PS)$-sequence $(z_n)$ is bounded in $E^\alpha(\Omega)$.
\end{lema}

{\bf Proof.} By (\ref{2.5}) we have for $z_n=(u_n, v_n)$
\[
\mathcal{J}(u_n, v_n) = \int\limits_\Omega\mathcal{A}^{\alpha/2}u_n\mathcal{A}^{s-\alpha/2}v_n dx -\frac{1}{p+1}\int\limits_\Omega v_n^{p+1} dx - \int\limits_\Omega F(u_n)dx\rightarrow c,
\]
and
\begin{equation}\label{2.7}
\vert\mathcal{J}'(u_n, v_n)(\varphi,\phi)\vert\leq\varepsilon_n\Vert(\varphi,\phi)\Vert_E\leq\varepsilon_n(\Vert \varphi\Vert_{\Theta^\alpha}+\Vert \phi\Vert_{\Theta^{2s-\alpha}}),
\end{equation}
where

\begin{eqnarray*}
\mathcal{J}'(u_n, v_n)(\varphi,\phi)&=&\int\limits_\Omega\mathcal{A}^{\alpha/2}u_n\mathcal{A}^{s-\alpha/2}\phi dx +\int\limits_\Omega\mathcal{A}^{s-\alpha/2}v_n\mathcal{A}^{\alpha/2}\varphi dx -\int\limits_\Omega v^p_n\phi dx\\
&&-\int\limits_\Omega f(u_n)\varphi dx.
\end{eqnarray*}

Choosing $(\varphi,\psi) = (u_n, v_n)\in \Theta^{\alpha}(\Omega)\times\Theta^{2s-\alpha}(\Omega)$ we get by (\ref{2.7})
\[
2\int\limits_\Omega\mathcal{A}^{\alpha/2}u_n\mathcal{A}^{s-\alpha/2}v_ndx -\int\limits_\Omega v_n^{p+1} dx - \int\limits_\Omega f(u_n)u_n dx\leq\varepsilon_n(\Vert u_n\Vert_{\Theta^\alpha}+\Vert v_n\Vert_{\Theta^{2s-\alpha}})
\]
and subtracting this from $2\mathcal{J}(u_n, v_n)$ we obtain, using assumption of Theorem \ref{teo1}
\[
\left(1-\frac{2}{p+1}\right)\int\limits_\Omega v_n^{p+1} dx +\left(1-\frac{2}{\theta}\right) \int\limits_\Omega f(u_n)u_n dx\leq C+\varepsilon_n(\Vert u_n\Vert_{\Theta^\alpha}+\Vert v_n\Vert_{\Theta^{2s-\alpha}})
\]
and thus
\begin{equation}\label{2.10}
\int\limits_\Omega v_n^{p+1} dx \leq C +\varepsilon_n(\Vert u_n\Vert_{\Theta^\alpha}+\Vert v_n\Vert_{\Theta^{2s-\alpha}}),
\end{equation}

\begin{equation}\label{2.11}
\int\limits_\Omega f(u_n)u_n dx\leq C +\varepsilon_n(\Vert u_n\Vert_{\Theta^\alpha}+\Vert v_n\Vert_{\Theta^{2s-\alpha}}).
\end{equation}
Choosing $(\varphi,\phi) = (0,\mathcal{A}^{\alpha-s}u_n)\in \Theta^{\alpha}(\Omega)\times\Theta^{2s-\alpha}(\Omega)$ in (\ref{2.7}) we get
\[
\int\limits_\Omega\vert\mathcal{A}^{\alpha/2}u_n\vert^2 dx \leq\int\limits_\Omega v^p_n\mathcal{A}^{\alpha-s}u_ndx +\varepsilon_n\Vert\mathcal{A}^{\alpha-s}u_n\Vert_{\Theta^{2s-\alpha}}
\]
and hence by H\"{o}lder inequality
\[
\Vert u_n\Vert^2_{\Theta^{\alpha}}=\Vert\mathcal{A}^{\alpha/2}u_n\Vert^2_{L^2}\leq \left(\int\limits_\Omega\vert v_n\vert^{p+1} dx\right)^{\frac{p}{p+1}}\left(\int\limits_\Omega \vert\mathcal{A}^{\alpha-s}u_n\vert^{p+1}dx\right)^{\frac{1}{p+1}} +\varepsilon_n\Vert u_n\Vert_{\Theta^{\alpha}}.
\]
Noting that
\[
\left(\int\limits_\Omega\vert\mathcal{A}^{\alpha-s}u_n\vert^{p+1} dx\right)^{\frac{1}{p+1}}\leq c\Vert \mathcal{A}^{\alpha-s}u_n\Vert_{\Theta^{2s-\alpha}} = c\Vert \mathcal{A}^{\alpha/2}u_n\Vert_{L^2} = c\Vert u_n\Vert_{\Theta^{\alpha}}
\]
we obtain, using (\ref{2.10})
\[
\Vert u_n\Vert^2_{\Theta^{\alpha}}\leq [C + \varepsilon_n(\Vert u_n\Vert_{\Theta^\alpha}+\Vert v_n\Vert_{\Theta^{2s-\alpha}})]^{p/(p+1)}\cdot c \Vert u_n\Vert_{\Theta^\alpha}+\varepsilon_n\Vert u_n\Vert_{\Theta^{\alpha}}
\]
and thus
\begin{equation}\label{2.12}
\Vert u_n\Vert_{\Theta^{\alpha}}\leq C + \varepsilon_n(\Vert u_n\Vert_{\Theta^\alpha}+\Vert v_n\Vert_{\Theta^{2s-\alpha}})^{p/(p+1)}.
\end{equation}
Similarly as above we note that $\mathcal{A}^{s-\alpha}v_n\in\Theta^\alpha(\Omega)$, and thus, choosing $(\varphi,\psi) =(\mathcal{A}^{s-\alpha}v_n, 0)\in \Theta^\alpha(\Omega)\times\Theta^{2s-\alpha}(\Omega)$ in (\ref{2.7}) we get
\begin{eqnarray*}
\int\limits_\Omega\vert \mathcal{A}^{s-\alpha/2}v_n\vert^2dx & \leq&\int\limits_\Omega f(u_n)\mathcal{A}^{s-\alpha}v_ndx + \varepsilon_n\Vert\mathcal{A}^{s-\alpha}v_n\Vert_{\Theta^\alpha}\\
&\leq &\Vert\mathcal{A}^{s-\alpha}v_n\Vert_\infty\int\limits_\Omega\vert f(u_n)\vert dx + \varepsilon_n\Vert v_n\Vert_{\Theta^{\alpha}}.
\end{eqnarray*}
Using that $\Vert\mathcal{A}^{s-\alpha}v_n\Vert_{\Theta^\alpha}= \Vert\mathcal{A}^{s-\alpha/2}v_n\Vert_{L^2} = \Vert v_n\Vert_{\Theta^{2s-\alpha}}$, and the fact that $\Theta^{\alpha}(\Omega)\subset C(\Omega)$ we then
obtain, using (\ref{2.11})
\begin{eqnarray*}\label{2.13}
\Vert v_n\Vert_{\Theta^\alpha}&\leq & c\int\limits_\Omega \vert f(u_n)\vert dx + \varepsilon_n \\
&= &\int\limits_{[\vert u_n\vert\leq s_0]}\max_{\vert t\vert\leq s_0}\vert f(t)\vert dx + \int\limits_{[\vert u_n\vert> s_0]} f(u_n)u_n dx + \varepsilon_n\\
&\leq & C + \varepsilon_n(\Vert u_n\Vert_{\Theta^{\alpha}} + \Vert v_n\Vert_{\Theta^{2s-\alpha}}).
\end{eqnarray*}
Joining (\ref{2.12}) and (\ref{2.13}) we finally get
\[
\Vert u_n\Vert_{\Theta^{\alpha}} + \Vert v_n\Vert_{\Theta^{2s-\alpha}}\leq C + 2\varepsilon_n(\Vert u_n\Vert_{\Theta^{\alpha}} + \Vert v_n\Vert_{\Theta^{2s-\alpha}}).
\]
Thus, $\Vert u_n\Vert_{\Theta^{\alpha}} + \Vert v_n\Vert_{\Theta^{2s-\alpha}}$ is bounded. \; \fim

With this it is now possible to complete the proof of the $(PS)$-condition: since $\Vert u_n\Vert_{\Theta^\alpha}$ is bounded, we find a weakly convergent subsequence $u_n\rightharpoonup u$ in $\Theta^{\alpha}(\Omega)$. Since the mappings $\mathcal{A}^{\alpha/2}: \Theta^{\alpha}(\Omega)\rightarrow L^2(\Omega)$ and $\mathcal{A}^{\alpha/2-s}: L^2(\Omega)\rightarrow \Theta^{2s-\alpha}(\Omega)$ are continuous isomorphisms, we get $\mathcal{A}^{\alpha/2}(u_n - u)\rightharpoonup 0$ in $L^2(\Omega)$ and $\mathcal{A}^{\alpha-s}(u_n - u) \rightharpoonup 0$ in $\Theta^{2s-\alpha}(\Omega)$. Since $\Theta^{2s-\alpha}(\Omega)\subset L^{p+1}(\Omega)$ compactly, we conclude that $\mathcal{A}^{\alpha-s}(u_n - u)\rightarrow 0$ strongly in $L^{p+1}(\Omega)$.

Similarly, we find a subsequence of $(v_n)$ which is weakly convergent in $\Theta^{2s-\alpha}(\Omega)$ and such that $v^p_n$ is strongly convergent in $L^{\frac{p+1}{p}}(\Omega)$.

Choosing $(\varphi,\phi) = (0,\mathcal{A}^{\alpha-s}(u_n - u))\in \Theta^{\alpha}(\Omega)\times\Theta^{2s-\alpha}(\Omega)$ in (\ref{2.7}) we thus conclude
\[
\int\limits_\Omega\mathcal{A}^{\alpha/2}u_n\mathcal{A}^{\alpha/2}(u_n - u)dx \leq\int\limits_\Omega v^p_n\mathcal{A}^{\alpha-s}(u_n - u)dx + \varepsilon_n\Vert\mathcal{A}^{\alpha-s}(u_n - u)\Vert_{\Theta^{2s-\alpha}}.
\]
By the above considerations, the right-hand-side converges to $0$, and thus
\[
\int\limits_\Omega\vert\mathcal{A}^{\alpha/2}u_n\vert^2 dx\rightarrow\int\limits_\Omega\vert\mathcal{A}^{\alpha/2}u\vert^2 dx.
\]
Thus, $u_n\rightarrow u$ strongly in $\Theta^{\alpha}(\Omega)$.

To obtain the strong convergence of $(v_n)$ in $\Theta^{2s-\alpha}(\Omega)$, one proceeds similarly: as
above, one finds a subsequence $(v_n)$ converging weakly in $\Theta^{2s-\alpha}(\Omega)$ to $v$, and then
$\mathcal{A}^{s-\alpha}v_n\rightharpoonup\mathcal{A}^{s-\alpha}v$ weakly in $\Theta^{\alpha}(\Omega)$ and $\mathcal{A}^{s-\alpha}v_n\rightarrow\mathcal{A}^{s-\alpha}v$ strongly in $C(\Omega)$. Choosing in (\ref{2.5}) $(\varphi,\phi) = (\mathcal{A}^{s-\alpha}(v_n - v), 0)$, we get
\[
\int\limits_\Omega \mathcal{A}^{s-\alpha/2}(v_n - v)\mathcal{A}^{s-\alpha/2}v_n dx \leq \int\limits_\Omega f(u_n)\mathcal{A}^{s-\alpha}(v_n - v)dx + \varepsilon_n(\Vert \mathcal{A}^{s-\alpha}(v_n - v)\Vert_{\Theta^{\alpha}}). 
\]
The first term on the right is estimated by
\[
\Vert \mathcal{A}^{s-\alpha}(v_n - v)\Vert_C\int\limits_\Omega\vert f(u_n)\vert dx\rightarrow 0,
\]
and thus one concludes again that
\[
\int\limits_\Omega\vert\mathcal{A}^{s-\alpha/2}v_n\vert^2 dx\rightarrow\int\limits_\Omega\vert\mathcal{A}^{s-\alpha/2}v\vert^2 dx
\]
and hence also $v_n\rightarrow v$ strongly in $\Theta^{2s-\alpha}(\Omega)$.

Thus, the conditions of Theorem \ref{teo2.4} are satisfied; hence, we find a positive
critical point $(u, v)$ for the functional $\mathcal{J}$, which yields a weak solution to system (\ref{1}).

\subsection{The case $p\leq 1$: Variational setting}

Suppose that $p\leq 1$, $n=2$ and $s\in\left(\frac{1}{2},1\right)$ or $n=3$ and $s\in\left(\frac{3}{4},1\right)$. Thus, $\Theta^{2s}(\Omega)$ is compactly embedded in $C(\Omega)$.

Let $\Omega$ be a smooth bounded open subset of $\R^n$ and $0 < s < 1$. In order to inspire our formulation, assume that the couple $(u,v)$ of nontrivial functions is roughly a solution of (\ref{1}). From the first equation, we have $v = \left(\mathcal{A}^{s} u\right)^{\frac1p}$. Plugging this equality into the second equation, we obtain

\begin{equation} \label{S3.4}
\left\{\begin{array}{rrll}
\mathcal{A}^{s} \left( \mathcal{A}^{s} u \right)^{\frac1p} &=& f(u) & {\rm in} \ \ \Omega\; \\
u &=& 0 & {\rm in} \ \ \partial\Omega
\end{array}\right..
\end{equation}

\n The basic idea in trying to solve (\ref{S3.4}) is considering the functional $\Phi : \Theta^{2s}(\Omega) \rightarrow \R$ defined by

\begin{equation}\label{func}
\Phi(u)=\frac{p}{p+1}\int\limits_{\Omega}\vert \mathcal{A}^{s}u\vert^{\frac{p+1}{p}}dx-\int\limits_{\Omega} F(u)dx\, .
\end{equation}

The Gateaux derivative of $\Phi$ at $u \in \Theta^{2s}(\Omega)$ in the direction $\varphi \in \Theta^{2s}(\Omega)$ is given by

\[
\Phi'(u)\varphi=\int\limits_{\Omega}\left|\mathcal{A}^{s}u\right|^{\frac{1}{p}-1} \mathcal{A}^{s}u \mathcal{A}^{s}\varphi dx - \int\limits_{\Omega}f(u) \varphi dx\,
\]
and thus, hence $f(u)\in C(\Omega)$ the problem

\[
\left\{\begin{array}{rrll}
\mathcal{A}^{s} v &=& f(u) & {\rm in} \ \ \Omega\; \\
v &=& 0 & {\rm in} \ \ \partial\Omega
\end{array}\right. 
\]
admits a unique nontrivial weak solution $v \in \Theta^s(\Omega)$. Then, one easily checks that $u$ is a weak solution of the problem
\[
\left\{\begin{array}{rrll}
\mathcal{A}^{s} u &=& v^p & {\rm in} \ \ \Omega\; \\
u &=& 0 & {\rm in} \ \ \partial\Omega
\end{array}\right. .
\]
In short, starting from a critical point $u \in \Theta^{2s}(\Omega)$ of $\Phi$, we have constructed a nontrivial weak solution $(u,v) \in \Theta^{2s}(\Omega) \times \Theta^s(\Omega)$ of the problem (\ref{1}).

\subsection{Existence of critical points}

From what we saw above, it suffices to show the existence of a nonzero critical point $u \in \Theta^{2s}(\Omega)$ of the functional $\Phi$.

In this case $p\leq 1$, $n=2$ and $s\in\left(\frac{1}{2},1\right)$ or $n=3$ and $s\in\left(\frac{3}{4},1\right)$. Thus, $\Theta^{2s}(\Omega)$ is compactly embedded in $C(\Omega)$. Then, the second term of the functional $\Phi$ is defined if $F$ is continuous, and no growth restriction on $F$ is necessary. Since $F$ is differentiable, the functional $\Phi$ is a well-defined $C^1$-functional on the space $\Theta^{2s}(\Omega)$. The proof consists in applying the classical mountain pass theorem of Ambrosetti and Rabinowitz in our variational setting. 

We now show that $\Phi$ has a local minimum in the origin. 

\begin{eqnarray*}
\Phi(u) &=& \frac{p}{p+1}\int\limits_{\Omega}\vert \mathcal{A}^{s}u\vert^{\frac{p+1}{p}}dx-\int\limits_{\Omega} F(u)dx\\
&\geq & \frac{p c}{p+1}\Vert u\Vert_{C}^{\frac{p+1}{p}} - o\left(\Vert u\Vert_{C}^{\frac{p+1}{p}}\right),
\end{eqnarray*}
so that the origin $u_0 = 0$ is a local minimum point. Next, let $u_1 = t \overline{u}$, where $t > 0$ and $\overline{u} \in \Theta^{2s}(\Omega)$ is a nonzero function. Then
\[
\Phi(u_1)\leq\frac{p t^{\frac{p+1}{p}}}{p+1} \int\limits_{\Omega} \vert \mathcal{A}^{s} \overline{u} \vert^{\frac{p+1}{p}}dx - t^{\theta} \Vert \overline{u}\Vert_{C}^{\theta}+d
\]
with $\theta>\frac{p+1}{p}$ (by assumption), and thus $\Phi(t\overline{u})\rightarrow-\infty$ as $t\rightarrow+\infty$.

Finally, we show that $\Phi$ fulfills the Palais-Smale condition (PS). Let $(u_k) \subset \Theta^{2s}(\Omega)$ be a (PS)-sequence, that is,
\[
\vert \Phi(u_k)\vert \leq C_0
\]
and

\[
\vert \Phi'(u_k) \varphi\vert \leq \varepsilon_k \Vert \varphi \Vert_{\Theta^{2s}}
\]
for all $\varphi \in \Theta^{2s}(\Omega)$, where $\varepsilon_k \rightarrow 0$ as $k \rightarrow +\infty$.

We have
\begin{eqnarray*}
C_0 + \varepsilon_k \Vert u_k \Vert_{\Theta^{2s}} &\geq & \vert \theta \Phi(u_k) - \Phi'(u_k) u_k\vert \\
&\geq &\left(\theta\frac{p}{p+1}-1\right) \int\limits_{\Omega}\vert \mathcal{A}^{s}u_k\vert^{\frac{p+1}{p}}dx-\theta\int\limits_{\Omega}F(u_k)dx+\int\limits_{\Omega}f(u_k)u_kdx\\
&\geq & \left(\theta\frac{p}{p+1}-1\right) \int\limits_{\Omega}\vert \mathcal{A}^{s}u_k\vert^{\frac{p+1}{p}}dx-C_0\geq \delta  \Vert u_k\Vert_{\Theta^{2s}}^{\frac{p+1}{p}}-C_0,
\end{eqnarray*}
and thus $(u_k)$ is bounded in $\Theta^{2s}(\Omega)$. Thanks to the compactness of the embedding $\Theta^{2s}(\Omega) \hookrightarrow C(\Omega)$, one easily checks that $(u_k)$ converges strongly in $\Theta^{2s}(\Omega)$. So, by the mountain pass theorem, we obtain a nonzero critical point $u \in \Theta^{2s}(\Omega)$. This ends the proof.

\section{Proof of theorem \ref{teo2}}

Note that if $n=4s$ we have $\Theta^{2s}(\Omega)$ is compactly embedded in $L^{r}(\Omega)$ for all $r>1$ and if $n>4s$ implies that $\Theta^{2s}(\Omega)$ is continuously embedded in $L^{\frac{2n}{n-4s}}(\Omega)$. It suffices to prove the result for $n > 4s$, since the ideas involved in its proof are fairly similar when $n = 4s$. 

\subsection{Variational setting}

Let $\Omega$ be a smooth bounded open subset of $\R^n$ and $0 < s < 1$. In order to inspire our formulation, assume that the couple $(u,v)$ of nonnegative functions is roughly a solution of (\ref{2}). From the first equation, we have $v = \left(\mathcal{A}^{s} u\right)^{\frac1p}$. Plugging this equality into the second equation, we obtain

\begin{equation} \label{S3.41}
\left\{\begin{array}{rrll}
\mathcal{A}^{s} \left( \mathcal{A}^{s} u \right)^{\frac1p} &=& u^q & {\rm in} \ \ \Omega\; \\
u &\geq& 0 & {\rm in} \ \ \Omega \\
u &=& 0 & {\rm in} \ \ \partial\Omega
\end{array}\right..
\end{equation}

\n The basic idea in trying to solve (\ref{S3.41}) is considering the functional $\Phi : \Theta^{2s}(\Omega) \rightarrow \R$ defined by

\begin{equation}\label{func}
\Phi(u)=\frac{p}{p+1}\int\limits_{\Omega}\vert \mathcal{A}^{s}u\vert^{\frac{p+1}{p}}dx-\frac{1}{q+1}\int\limits_{\Omega} (u^+)^{q+1}dx\, .
\end{equation}

The Gateaux derivative of $\Phi$ at $u \in \Theta^{2s}(\Omega)$ in the direction $\varphi \in \Theta^{2s}(\Omega)$ is given by

\[
\Phi'(u)\varphi=\int\limits_{\Omega}\left|\mathcal{A}^{s}u\right|^{\frac{1}{p}-1} \mathcal{A}^{s}u \mathcal{A}^{s}\varphi dx - \int\limits_{\Omega}(u^+)^q \varphi dx\, .
\]
In this case, $\Theta^{2s}(\Omega)$ is continuously embedded in $L^{\frac{2n}{n-4s}}(\Omega)$. Thus, if $0<q\leq\frac{n+2s}{n-4s}$ we have $u^q\in L^{\frac{2n}{n+2s}}(\Omega)$. Therefore the problem
\begin{equation}\label{p}
\left\{\begin{array}{rrll}
\mathcal{A}^{s} v &=& (u^+)^q & {\rm in} \ \ \Omega\; \\
v &=& 0 & {\rm in} \ \ \partial\Omega
\end{array}\right.
\end{equation}
admits a unique nonnegative weak solution $v \in \Theta^s(\Omega)$. Now, if $\frac{n+2s}{n-4s}<q<\frac{n+4s}{n-4s}$ then $u^q\in \Theta^{r-2s}(\Omega)$, where $0<r:=\frac{n+4s-(n-4s)q}{2}<s$. Therefore the problem (\ref{p}) admits a unique nonnegative weak solution $v \in \Theta^r(\Omega)$. 

Then, one easily checks that $u$ is a weak solution of the problem
\[
\left\{\begin{array}{rrll}
\mathcal{A}^{s} u &=& v^p & {\rm in} \ \ \Omega\; \\
u &=& 0 & {\rm in} \ \ \partial\Omega
\end{array}\right. .
\]
In short, starting from a critical point $u \in \Theta^{2s}(\Omega)$ of $\Phi$, we have constructed a nonnegative weak solution 
\[
(u,v)\in\left\{\begin{array}{rrll}
\Theta^{2s}(\Omega) \times \Theta^s(\Omega)  & {\rm if} \ \ 0<q\leq\frac{n+2s}{n-4s}\; \\
\Theta^{2s}(\Omega) \times \Theta^r(\Omega)  & {\rm if} \ \ \frac{n+2s}{n-4s}<q<\frac{n+4s}{n-4s}
\end{array}\right.
\]
of the problem (\ref{2}).

\subsection{The existence part in case $pq < 1$}
We apply the direct method to the functional $\Phi$ on $\Theta^{2s}(\Omega)$.

In order to show the coercivity of $\Phi$, note that $q + 1 < \frac{p+1}{p}$ because $pq < 1$. Hence $q<\frac{n+4s}{n-4s}$ the embedding $\Theta^{2s}(\Omega) \hookrightarrow L^{q+1}(\Omega)$ is continuous. So, for $p\leq 1$ there exist constants $C_1, C_2 > 0$ such that
\begin{eqnarray*}
\Phi(u) &=& \frac{p}{p+1}\int\limits_{\Omega}\vert \mathcal{A}^{s}u\vert^{\frac{p+1}{p}}dx - \frac{1}{q+1}\int\limits_{\Omega}\vert u\vert ^{q+1}dx\\
&\geq& \frac{p C_1}{p+1}\Vert u\Vert_{\Theta^{2s}}^{\frac{p+1}{p}}-\frac{C_2}{q+1}\Vert u\Vert_{\Theta^{2s}}^{q+1}\\
&=& \Vert u\Vert^{\frac{p+1}{p}}_{\Theta^{2s}}\left(\frac{p C_1}{p+1}-\frac{C_2}{(q+1)\Vert u\Vert^{\frac{p+1}{p}-(q+1)}_{\Theta^{2s}}}\right)
\end{eqnarray*}
for all $u\in \Theta^{2s}(\Omega)$. Therefore, $\Phi$ is lower bounded and coercive, that is, $\Phi(u) \rightarrow +\infty$ as $\Vert u\Vert_{\Theta^{2s}} \rightarrow +\infty$.

Let $(u_k) \subset \Theta^{2s}(\Omega)$ be a minimizing sequence of $\Phi$. It is clear that $(u_k)$ is bounded in $\Theta^{2s}(\Omega)$, since $\Phi$ is coercive. So, module a subsequence, we have $u_k \rightharpoonup u_0$ in $\Theta^{2s}(\Omega)$. Since $\Theta^{2s}(\Omega)$ is compactly embedded in $L^{q+1}(\Omega)$, we have $u_k \rightarrow u_0$ in $L^{q+1}(\Omega)$. Here, we again use the fact that $q + 1 < \frac{p+1}{p}$. Thus,
\begin{eqnarray*}
\lim_{n\rightarrow\infty}\inf \Phi(u_k) &=& \lim_{k \rightarrow \infty} \inf \frac{p}{p+1}\Vert \mathcal{A}^{s} u_k \Vert^{\frac{p+1}{p}}_{L^{\frac{p+1}{p}}}-\frac{1}{q+1}\Vert u_{0}\Vert_{L^{q+1}}^{q+1}\\
&\geq &\frac{p}{p+1}\Vert \mathcal{A}^{s}u_{0}\Vert^{\frac{p+1}{p}}_{L^{\frac{p+1}{p}}}-\frac{1}{q+1}\Vert u_{0}\Vert_{L^{q+1}}^{q+1} = \Phi(u_{0})\, ,
\end{eqnarray*}
so that $u_0$ minimizers $\Phi$ on $\Theta^{2s}(\Omega)$. We just need to guarantee that $u_0$ is nonzero. But, this fact is clearly true since $\Phi(\varepsilon u_1) < 0$ for any nonzero nonnegative function $u_1 \in \Theta^{2s}(\Omega)$ and $\varepsilon > 0$ small enough, that is,

\[
\Phi(\varepsilon u_1) = \frac{p\varepsilon^{\frac{p+1}{p}}}{p+1}\int\limits_{\Omega}\vert \mathcal{A}^{s} u_1 \vert^{\frac{p+1}{p}}dx - \frac{\varepsilon^{q+1}}{q+1}\int\limits_{\Omega}\vert u_1 \vert ^{q+1}dx < 0
\]
for $\varepsilon>0$ small enough. This ends the proof of existence.\; \fim

\subsection{The uniqueness part in case $pq < 1$}
The main tools in the proof of uniqueness are the strong maximum principle and a Hopf's lemma adapted to fractional operators.

Let $(u_{1},v_{1}),(u_{2},v_{2})$ be two positive solutions of (\ref{2}). Define
\[
S=\{s\in(0,1] : u_{1} - tu_{2},\ v_{1} - tv_{2}\geq 0\ \text{ in } \overline{\Omega}\ \text{ for all }\ t \in [0,s]\}\, .
\]
From the strong maximum principle and Hopf's lemma (see \cite{T}), it follows that $S$ is no empty.

Let $s_{\ast}=\sup S$ and assume that $s_{\ast} < 1$.

Clearly,

\begin{equation}\label{eq1}
u_{1}-s_{\ast}u_{2},\ v_{1}-s_{\ast}v_{2}\geq 0\ \text{ in }\ \overline{\Omega}\, .
\end{equation}
By (\ref{eq1}) and the integral representation in terms of the Green function $G_{\Omega}$ of $\mathcal{A}^s$ (see \cite{Green, ros176}), we have

\begin{eqnarray*}
u_{1}(x)&=&\int\limits_{\Omega}G_{\Omega}(x,y)v_{1}^{p}(y)dy\geq\int\limits_{\Omega}G_{\Omega}(x,y)s_{\ast}^p v_{2}^{p}(y)dy\\
&=& s_{\ast}^{p}\int\limits_{\Omega}G_{\Omega}(x,y)v_{2}^{p}(y)dy=s_{\ast}^{p}u_{2}(x)
\end{eqnarray*}
for all $x \in \overline{\Omega}$. In a similar way, one gets $v_1 \geq s_{\ast}^{q} v_2$ in $\overline{\Omega}$.

\n Using the assumption $pq < 1$ and the fact that $s_{\ast} < 1$, we derive

\begin{equation}
\left\{\begin{array}{llll}
\mathcal{A}^s(u_{1} - s_{\ast}u_{2}) = v_{1}^{p}-s_{\ast}v_{2}^{p} \geq (s_{\ast}^{pq}-s_{\ast})v_{2}^{p} > 0  \\
\mathcal{A}^s(v_{1} - s_{\ast}v_{2}) = u_{1}^{q}-s_{\ast}u_{2}^{q} \geq (s_{\ast}^{pq}-s_{\ast})u_{2}^{q} > 0 \ \
\end{array}\right.\ {\rm in}\ \Omega
\end{equation}
So, by the strong maximum principle, one has $u_{1}-s_{\ast}u_{2}, v_{1}-s_{\ast}v_{2} > 0$ in $\Omega$. Then, by Hopf's lemma, we have $\frac{\partial}{\partial\nu}(u_{1}-s_{\ast}u_{2}),\frac{\partial}{\partial\nu}(v_{1}-s_{\ast}v_{2})<0$ on $\partial \Omega$, where $\nu$ is the unit outer normal in $\R^n$ to $\partial\Omega$, so that $u_{1} - (s_{\ast}+\varepsilon)u_{2}, v_{1} - (s_{\ast}+\varepsilon)v_{2} > 0$ in $\Omega$ for $\varepsilon > 0$ small enough, contradicting the definition of $s_{\ast}$. Therefore, $s_{\ast} \geq 1$ and, by (\ref{eq1}), $u_{1} - u_{2}, v_{1} - v_{2} \geq 0$ in $\overline{\Omega}$. A similar reasoning also produces $u_{2} - u_{1}, v_{2} - v_{1} \geq 0$ in $\overline{\Omega}$. This ends the proof of uniqueness.\; \fim

\subsection{Existence of critical points in case $pq>1$}

From what we saw, it suffices to show the existence of a nonzero critical point $u \in \Theta^{2s}(\Omega)$ of the functional $\Phi$.

Again the proof consists in applying the classical mountain pass theorem of Ambrosetti and Rabinowitz in our variational setting. We first assert that $\Phi$ has a local minimum in the origin. 

Note that $p\leq 1$ and $q + 1 > \frac{p+1}{p}$ because $pq > 1$. Hence $q<\frac{n+4s}{n-4s}$ the embedding $\Theta^{2s}(\Omega) \hookrightarrow L^{q+1}(\Omega)$ is compact. Consider the set $\Gamma :=\left\{u\in \Theta^{2s}(\Omega) : \Vert u\Vert_{\Theta^{2s}}=\rho\right\}$. Then, on $\Gamma$, we have

\begin{eqnarray*}
\Phi(u) &=& \frac{p}{p+1}\int\limits_{\Omega}\vert \mathcal{A}^{s}u\vert^{\frac{p+1}{p}}dx-\frac{1}{q+1}\int\limits_{\Omega}\vert u\vert^{q+1}dx\\
&\geq & \frac{p C_1}{p+1}\Vert u\Vert_{\Theta^{2s}}^{\frac{p+1}{p}} - \frac{C_2}{q+1}\Vert u\Vert_{\Theta^{2s}}^{q+1} =\rho^{\frac{p+1}{p}}\left( \frac{p C_1}{p+1}-\frac{C_2}{q+1}\rho^{q+1-\frac{p+1}{p}}\right) \\
&>& 0 = \Phi(0)
\end{eqnarray*}
for fixed $\rho>0$ small enough, so that the origin $u_0 = 0$ is a local minimum point. In particular,  $\inf_{\Gamma} \Phi > 0 = \Phi(u_0)$.

Note that $\Gamma$ is a closed subset of $\Theta^{2s}(\Omega)$ and decomposes $\Theta^{2s}(\Omega)$ into two connected components, namely $\left\{u\in \Theta^{2s}(\Omega) : \Vert u\Vert_{\Theta^{2s}} < \rho \right\}$ and $\left\{u\in \Theta^{2s}(\Omega) : \Vert u\Vert_{\Theta^{2s}} > \rho\right\}$.

Let $u_1 = t \overline{u}$, where $t > 0$ and $\overline{u} \in \Theta^{2s}(\Omega)$ is a nonzero nonnegative function. Since $pq > 1$, we can choose $t$ sufficiently large so that
\[
\Phi(u_1)=\frac{p t^{\frac{p+1}{p}}}{p+1} \int\limits_{\Omega} \vert \mathcal{A}^{s} \overline{u} \vert^{\frac{p+1}{p}}dx - \frac{t^{q+1}}{q+1} \int\limits_{\Omega} (\overline{u}^+)^{q+1}dx < 0\, .
\]
It is clear that $u_1 \in \left\{u\in \Theta^{2s}(\Omega): \Vert u\Vert_{\Theta^{2s}} > \rho\right\}$. Moreover, $\inf_{\Gamma} \Phi > \max\{\Phi(u_0), \Phi(u_1)\}$, so that the mountain pass geometry is satisfied.

Finally, we show that $\Phi$ fulfills the Palais-Smale condition (PS). Let $(u_k) \subset \Theta^{2s}(\Omega)$ be a (PS)-sequence, that is,
\[
\vert \Phi(u_k)\vert \leq C_0
\]
and

\[
\vert \Phi'(u_k) \varphi\vert \leq \varepsilon_k \Vert \varphi \Vert_{\Theta^{2s}}
\]
for all $\varphi \in \Theta^{2s}(\Omega)$, where $\varepsilon_k \rightarrow 0$ as $k \rightarrow +\infty$.

From these two inequalities and the assumption $pq > 1$, we deduce that
\begin{eqnarray*}
C_0 + \varepsilon_k \Vert u_k \Vert_{\Theta^{2s}} &\geq & \vert (q+1) \Phi(u_k) - \Phi'(u_k) u_k\vert \\
&\geq &\left(\frac{p(q+1)}{p+1}-1\right) \int\limits_{\Omega}\vert \mathcal{A}^{s}u_k\vert^{\frac{p+1}{p}}dx\\
&\geq & C \Vert u_k\Vert_{\Theta^{2s}}^{\frac{p+1}{p}}
\end{eqnarray*}
and thus $(u_k)$ is bounded in $\Theta^{2s}(\Omega)$. Thanks to the compactness of the embedding $\Theta^{2s}(\Omega) \hookrightarrow L^{q+1}(\Omega)$, one easily checks that $(u_k)$ converges strongly in $\Theta^{2s}(\Omega)$. So, by the mountain pass theorem, we obtain a nonzero critical point $u \in \Theta^{2s}(\Omega)$. This ends the proof.

\section{The regularity part}

In this section, we establish regularity property of the weak solutions of system (\ref{1}) based on the results obtained in \cite{CS} and \cite{T}. Also we establish the Brezis-Kato type result and derive the regularity of solutions to (\ref{2}).

\begin{propo}
Let $(u,v)$ be a weak solution of the problem (\ref{1}). In the hypothesis of Theorem 1, we have $(u, v) \in L^{\infty}(\Omega) \times L^{\infty}(\Omega)$ and, moreover, $(u,v) \in C^{1,\beta}(\overline{\Omega}) \times C^{\beta}(\overline{\Omega})$ for some $\beta \in (0, 1)$. Now if $f$ be a $C^1$ function such that $f(0)=0$ we have $(u,v) \in C^{1,\beta}(\overline{\Omega}) \times C^{1,\beta}(\overline{\Omega})$ for some $\beta \in (0, 1)$.
\end{propo}
\n {\bf Proof.} In the case $p>1$ we find solution $(u,v)\in\Theta^{\alpha}(\Omega) \times \Theta^{2s-\alpha}(\Omega)$. By choosing $\alpha$ (see Lemma \ref{lema3.1}), and by Sobolev imbedding theorem (see \cite{DD}) we have $u\in L^{\infty}(\Omega)$. Then $f(u)\in L^{\infty}(\Omega)$. Thus, by regularity result (see \cite{CS} or Proposition 3.1 of \cite{T}) we have $v\in C^{\gamma}(\overline{\Omega})$ for some $\gamma \in (0, 1)$. Hence $\gamma +2s>1$ and $v^p\in C^{\gamma}(\overline{\Omega})$ again by Proposition 3.1 of \cite{T} we have $u\in C^{1,\gamma+2s-1}(\overline{\Omega})$. Therefore $(u,v) \in C^{1,\beta}(\overline{\Omega}) \times C^{\beta}(\overline{\Omega})$ for some $\beta \in (0, 1)$.

Now if $f$ be a $C^1$ function such that $f(0)=0$ analogously we have $(u,v) \in C^{1,\gamma+2s-1}(\overline{\Omega}) \times C^{\gamma}(\overline{\Omega})$ for some $\gamma \in (0, 1)$. Then $f(u)\in C^{2s}(\overline{\Omega})$. Hence $2s +2s>1$ again by Proposition 3.1 of \cite{T} we have $v\in C^{1,2s+2s-1}(\overline{\Omega})$.  Therefore $(u,v) \in C^{1,\beta}(\overline{\Omega}) \times C^{1,\beta}(\overline{\Omega})$ for some $\beta \in (0, 1)$.

In the case $p\leq 1$ we find solution $(u,v)\in\Theta^{2s}(\Omega) \times \Theta^{s}(\Omega)$. From Sobolev imbedding theorem (see \cite{DD}) we have $u\in L^{\infty}(\Omega)$. Analogous to the previous case, we have the result.\ \fim\\

Next we prove the $L^\infty$ estimate of Brezis-Kato type.

\begin{propo}\label{prop}
Let $(u,v)$ be a weak solution of the problem (\ref{2}). In the hypothesis of Theorem 2, we have $(u, v) \in L^{\infty}(\Omega) \times L^{\infty}(\Omega)$ and, moreover, $(u,v) \in C^{1,\beta}(\overline{\Omega}) \times C^{1,\beta}(\overline{\Omega})$ for some $\beta \in (0, 1)$.
\end{propo}
\n {\bf Proof.} It suffices to prove the result for $n > 4s$, since the ideas involved in its proof are fairly similar when $n = 4s$.

In this case we find solution 
\[
(u,v)\in\left\{\begin{array}{rrll}
\Theta^{2s}(\Omega) \times \Theta^s(\Omega)  & {\rm if} \ \ 0<q\leq\frac{n+2s}{n-4s}\; \\
\Theta^{2s}(\Omega) \times \Theta^r(\Omega)  & {\rm if} \ \ \frac{n+2s}{n-4s}<q<\frac{n+4s}{n-4s}
\end{array}\right. ,
\]
where $0<r=\frac{n+4s-(n-4s)q}{2}<s$.

We analyze separately two different cases depending on the values of $q$. Note that $0< p \leq 1$.

For $q>1$, we rewrite the problem (\ref{2}) as follows
\begin{equation}
\left\{
\begin{array}{llll}
\mathcal{A}^{s}u = a(x)v^{\frac{p}{2}} & {\rm in} \ \ \Omega\\
\mathcal{A}^{s}v = b(x)u & {\rm in} \ \ \Omega\\
u= v=0 & {\rm in} \ \ \R^n\setminus\Omega
\end{array}
\right.,
\end{equation}
where $a(x) = v(x)^{\frac{p}{2}}$ and $b(x) = u(x)^{q - 1}$. Note that $p+1<\frac{2n}{n-2s}$ and $p+1\leq 2<\frac{2n}{n-2r}$. By Sobolev embedding, $\Theta^{s}(\Omega)\hookrightarrow L^{p+1}(\Omega)$ and $\Theta^{r}(\Omega)\hookrightarrow L^{p+1}(\Omega)$, so that $a \in L^{\frac{2(p+1)}{p}}(\Omega)$. Thus, for each fixed $\varepsilon > 0$, we can construct functions $q_{\varepsilon} \in L^{\frac{2(p+1)}{p}}(\Omega)$, $f_{\varepsilon} \in L^{\infty}(\Omega)$ and a constant $K_{\varepsilon} > 0$ such that

\[
a(x) v(x)^{\frac{p}{2}} = q_{\varepsilon}(x)v(x)^{\frac{p}{2}} + f_{\varepsilon}(x)
\]
and

\[
\Vert q_{\varepsilon} \Vert_{L^{\frac{2(p+1)}{p}}} < \varepsilon,\ \ \Vert f_{\varepsilon} \Vert_{L^{\infty}} < K_{\varepsilon}\, .
\]
In fact, consider the set

\[
\Omega_k = \{x \in\Omega: \vert a(x)\vert < k\}\, ,
\]
where $k$ is chosen such that

\[
\int\limits_{\Omega_k^c} \vert a(x)\vert^{\frac{2(p+1)}{p}}dx < \frac{1}{2}\varepsilon^{\frac{2(p+1)}{p}}\, .
\]
This condition is clearly satisfied for $k = k_\varepsilon$ large enough.

We now write

\begin{equation}
q_{\varepsilon}(x)=\left\{
\begin{array}{llll}
\frac{1}{m}a(x) & {\rm for} \ \ x \in \Omega_{k_\varepsilon}\\
a(x) & {\rm for} \ \ x \in \Omega_{k_\varepsilon}^c\\
\end{array}
\right.
\end{equation}
and

\[
f_{\varepsilon}(x) = \left(a(x) - q_{\varepsilon}(x)\right)v(x)^{\frac{p}{2}}\, .
\]
Then,

\begin{eqnarray*}
\int\limits_{\Omega}\vert q_{\varepsilon}(x)\vert^{\frac{2(p+1)}{p}}dx&=&\int\limits_{\Omega_{k_\varepsilon}}\vert q_{\varepsilon}(x)\vert^{\frac{2(p+1)}{p}}dx+\int\limits_{\Omega_{k_\varepsilon}^c}\vert q_{\varepsilon}(x)\vert^{\frac{2(p+1)}{p}}dx\\
&=&\left(\frac{1}{m}\right)^{\frac{2(p+1)}{p}}\int\limits_{\Omega_{k_\varepsilon}}\vert a(x)\vert^{\frac{2(p+1)}{p}}dx+\int\limits_{\Omega_{k_\varepsilon}^c}\vert a(x)\vert^{\frac{2(p+1)}{p}}dx\\
&<&\left(\frac{1}{m}\right)^{\frac{2(p+1)}{p}}\int\limits_{\Omega_{k_\varepsilon}}\vert a(x)\vert^{\frac{2(p+1)}{p}}dx+\frac{1}{2}\varepsilon^{\frac{2(p+1)}{p}}\, .
\end{eqnarray*}
So, for $m = m_\varepsilon > \left(\frac{2^{\frac{p}{2(p+1)}}}{\varepsilon}\right) \Vert a\Vert_{L^{\frac{2(p+1)}{p}}}$, we get

\[
\Vert q_{\varepsilon}\Vert_{L^{\frac{2(p+1)}{p}}} < \varepsilon\, .
\]
Note also that $f_{\varepsilon}(x) = 0$ for all $x \in \Omega_{k_\varepsilon}^c$ and, for this choice of $m$,

\[
f_{\varepsilon}(x) = \left( 1-\frac{1}{m_\varepsilon}\right) a(x)^2 \leq \left( 1-\frac{1}{m_\varepsilon}\right) k_\varepsilon^2
\]
for all $x \in \Omega_{k_\varepsilon}$. Therefore,

\[
\Vert f_{\varepsilon} \Vert_{L^{\infty}} \leq \left( 1-\frac{1}{m_\varepsilon}\right) k_\varepsilon^2 := K_\varepsilon\, .
\]
On the other hand, we have

\[
v(x) = \mathcal{A}^{-s}(bu)(x)\, ,
\]
where $b\in L^{\frac{q+1}{q-1}}(\Omega)$. Hence,

\[
u(x) = \mathcal{A}^{-s}\left[q_{\varepsilon}(x)(\mathcal{A}^{-s}(bu)(x))^{\frac{p}{2}}\right] + \mathcal{A}^{-s}f_{\varepsilon}(x)\, .
\]

By Lemma 2.1 of \cite{choi}, the claims $(ii)$ and $(iv)$ below follow readily and, by using H\"{o}lder's inequality, we also get the claims $(i)$ and $(iii)$. Precisely, for fixed $\gamma > 1$, we have:

 \begin{itemize}
    \item[(i)] The map $w \rightarrow b(x)w$ is bounded from $L^\gamma(\Omega)$ to $L^\beta(\Omega)$ for
\[
\frac{1}{\beta} = \frac{q-1}{q+1} + \frac{1}{\gamma};
\]
    \item[(ii)] For $\theta$ given by
\[
2s=n\left(\frac{1}{\beta}-\frac{2}{p \theta}\right),
\]
there exists a constant $C > 0$, depending on $\beta$ and $\theta$, such that

\[
\Vert (\mathcal{A}^{s}w)^{\frac{p}{2}}\Vert_{L^\theta}\leq C\Vert w\Vert_{L^\beta}^{\frac{p}{2}}
\]
for all $w \in L^\beta(\Omega)$;

    \item[(iii)] The map $ w \rightarrow q_\varepsilon(x)w$ is bounded from $L^\theta(\Omega)$ to $L^\eta(\Omega)$ with norm given by $\Vert q_\varepsilon \Vert_{L^{\frac{2(p+1)}{p}}}$, where $\theta \geq 1$ and $\eta$ satisfies
\[
\frac{1}{\eta}=\frac{p}{2(p+1)}+\frac{1}{\theta};
\]
    \item[(iv)] For $\delta$ given by

\[
2s = n\left(\frac{1}{\eta} - \frac{1}{\delta}\right),
\]
the map $ w\rightarrow \mathcal{A}^{-s}w$ is bounded from $L^\eta(\Omega)$ to $L^\delta(\Omega)$.
    \end{itemize}

Joining $(i)$, $(ii)$, $(iii)$ and $(iv)$, one easily checks that $\gamma < \delta$ and, in addition,

\begin{eqnarray*}
\Vert u\Vert_{L^\delta} &\leq & \Vert \mathcal{A}^{-s}\left[q_\varepsilon(x)\left(\mathcal{A}^{-s}(bu)\right)^{\frac{p}{2}}\right]\Vert_{L^\delta}+\Vert \mathcal{A}^{-s}f_\varepsilon\Vert_{L^\delta}\\
&\leq & C \left(\Vert q_\varepsilon\Vert_{L^{\frac{2(p+1)}{p}}}\Vert u\Vert_{L^\delta}^{\frac{p}{2}} + \Vert f_\varepsilon\Vert_{L^\delta}\right).
\end{eqnarray*}
Using now the fact that $p \leq 1$, $\Vert q_\varepsilon \Vert_{L^{\frac{2(p+1)}{p}}} < \varepsilon$ and $f_\varepsilon \in L^\infty(\Omega)$, we deduce that $\Vert u \Vert_{L^\delta} \leq C$ for some constant $C > 0$ independent of $u$. Proceeding inductively, we get $u \in L^\delta(\Omega)$ for all $\delta \geq 1$. So, Lemma 2.1 of \cite{choi} implies that $v \in L^{\infty}(\Omega)$. From this, and using Lemma 2.1 of \cite{choi} again, we deduce that $u \in L^{\infty}(\Omega)$.

Then $u^q,v^p\in L^{\infty}(\Omega)$. Thus, by regularity result (see \cite{CS} or Proposition 3.1 of \cite{T}) we have $v,u\in C^{2s}(\overline{\Omega})$. Hence it holds that $u^q,v^p \in C^{2s}(\overline{\Omega})$. Again, we can apply regularity result to deduce that $u\in C^{4s}(\overline{\Omega})$. Iteratively, we can raise the regularity so that $(u,v) \in C^{1,\beta}(\overline{\Omega}) \times C^{1,\beta}(\overline{\Omega})$ for some $\beta \in (0, 1)$.

The other case is treated in a similar way by writing $b(x) = u(x)^{\frac{q}{2}}$ if $q \leq 1$.\ \fim

\end{document}